\documentclass[11pt]{amsart}
\usepackage{amssymb}
\usepackage{graphicx}

\renewcommand{\a}{\alpha}

\newcommand{\g}{\gamma}
\newcommand{\e}{\varepsilon}

\newcommand{\SA}{{\mathcal{A}}}

\newcommand{\SC}{{\mathcal{C}}}

\newcommand{\SF}{{\mathcal{F}}}
\newcommand{\SG}{{\mathcal{G}}}

\newcommand{\SJ}{{\mathcal{J}}}

\newcommand{\SL}{{\mathcal{L}}}

\newcommand{\SO}{{\operatorname{SO}}}
\newcommand{\SP}{{\mathcal{P}}}

\renewcommand{\SS}{{\mathcal{S}}}

\newcommand{\C}{\mathbb{C}}

\newcommand{\R}{\mathbb{R}}

\newcommand{\Hom}{\operatorname{Hom}}
\newcommand{\length}{\operatorname{length}}
\newcommand{\id}{\textbf{\textit{I}}}

\newcommand{\End}{\operatorname{End}}

\newcommand{\GL}{\operatorname{GL}}
\newcommand{\Id}{\operatorname{Id}}
\newcommand{\Diff}{\operatorname{Diff}}
\newcommand{\grad}{\operatorname{grad}}

\renewcommand{\Re}{\operatorname{Re}}

\newcommand{\so}{{\mathfrak s\mathfrak o}}
\newcommand{\gl}{{\mathfrak g\mathfrak l}}

\newtheorem{proposition}{Proposition}[section]
\newtheorem{theorem}[proposition]{Theorem}

\newtheorem{lemma}[proposition]{Lemma}

\newtheorem{corollary}[proposition]{Corollary}
\newtheorem{remark}[proposition]{Remark}

\begin{document}

\title{Geometric structures on loop and path spaces}

\subjclass{Primary: 58B20. Secondary: 53D35, 55P35.}
\date{November, 2007}

\keywords{Loop space, symplectic structures.}
\thanks{First author partially supported through grant MEC
(Spain) MTM2007-63582}

\author{Vicente Mu\~noz}
\address{Departamento de Matem\'aticas \\ Consejo Superior de Investigaciones Cient\'{\i}ficas
\\ 28006 Madrid \\ Spain}
\address{Facultad de Matem\'aticas\\ Universidad Complutense de Madrid\\
28040 Madrid\\ Spain} \email{vicente.munoz@imaff.cfmac.csic.es}

\author{Francisco Presas}
\address{Departamento de Matem\'aticas \\
Consejo Superior de Investigaciones Cient\'{\i}ficas
\\ 28006 Madrid \\ Spain}
\address{Facultad de Matem\'aticas\\ Universidad Complutense de Madrid\\ 28040 Madrid\\ Spain}
\email{fpresas@imaff.cfmac.csic.es}


\renewcommand{\theenumi}{\roman{enumi}}

\begin{abstract}
Is is known that the loop space associated to a Riemannian manifold
admits a quasi-symplectic structure. This article shows that this
structure is not likely to recover the underlying Riemannian metric
by proving a result that is a strong indication of the ``almost''
independence of the quasi-symplectic structure with respect to the
metric. Finally conditions to have contact structures on these
spaces are studied.
\end{abstract}

\maketitle

\section{Introduction} \label{sect:introduction}
The loop space $\SL(M)$ of a manifold $M$ comes equipped with a
natural section of its tangent bundle defined as
\begin{eqnarray*}
\a: \SL(M) & \to & T\SL(M) \\
\gamma & \to &\gamma'. \
\end{eqnarray*}
Whenever we fix a Riemannian metric $g$ on $M$ we can define an
associated metric on the space of loops as
$$ (g_{\SL})_{\g}(X, Y)= \int_0^1 g(X(t), Y(t))dt, $$
where $X,Y\in T_\gamma \SL(M) \simeq \Gamma(S^1, \gamma^* TM)$ are
two tangent vectors. This metric gives us an isomorphism between
$T\SL(M)$ and $T^*\SL(M)$. Therefore $\a$ and $g_{\SL}$ allow us to
define a $1$-form
\begin{equation}
\mu(X)= \frac12 \int_0^1 g(X(t), \gamma'(t))dt, \label{eq:one_form}
\end{equation}
whose exterior differential $\omega= d\mu$ happens to be
quasi-symplectic. This means that the kernel of the form is finite
dimensional, specifically
$$ \ker (\omega)_{\gamma} = \{ X\in \Gamma(S^1, \gamma^*TM)\ ; \ \nabla_{\gamma'} X=0 \}. $$
(We assume throughout this article that the spaces considered are in
the $C^{\infty}$ category and have a natural Fr\'echet structure,
unless something else is declared.)

This quasi-symplectic structure can be enriched in many cases. This
is well known in the case of loop groups (i.e., $M$ is a Lie group). In this case
it is possible to define an integrable complex structure making a
finite codimensional closed manifold of a loop group into a
K\"ahler manifold.

Extending our space to the path space defined as
$$ \SP(M) = \{\gamma:[0,1] \to M \},$$
we still have the canonical section of the tangent bundle given by
\begin{eqnarray*}
\a: \SP(M) & \to & T\SP(M) \\
\gamma & \to &\gamma'. \
\end{eqnarray*}
As in the loop space we will easily check that equation
(\ref{eq:one_form}) is a $1$-form whose differential is symplectic.

\begin{proposition} \label{thm:str_path}
The $2$-form $\omega= d\mu$ in $\SP(M)$ induces a symplectic
structure. Moreover $(\SL(M), \omega)$ is a closed
quasi-symplectic submanifold of $\SP(M)$.
\end{proposition}

It could be thought that the symplectic structure makes life easier
in comparison to the quasi-symplectic one. At least, in terms of the
stability of the structure this is not the case. In particular, we
prove

\begin{theorem} \label{thm:main_unique}
Fix a smooth closed manifold $M$ of even dimension. Denote as $\omega_g$
the quasi-symplectic form on $\SL(M)$ associated to a Riemannian
metric $g$ on $M$. Then given two Riemannian metrics $g_0$ and
$g_1$ on $M$, there exists a smooth isotopy $\phi_{\epsilon}$ on
$\SL(M)$ such that $(\phi_{\epsilon})_{*}\omega_{g_0}$ is
$\epsilon$-close to $\omega_{g_1}$ in $L^2$-norm.
\end{theorem}

The proof of this result cannot be generalized to the case of
$\SP(M)$. Also, the theorem cannot be improved to obtain an
isotopy that matches $\omega_{g_0}$ and $\omega_{g_1}$ on the
nose. Certainly, allowing $\epsilon$ go to zero makes the norm of
the isotopy go to infinity and so discontinuities are developed.
So the result can be understood as a sort of ``approximate
uniqueness'' of the quasi-symplectic structure.

The proof of Theorem \ref{thm:main_unique} is based on an
adaptation of Moser's trick to this setting. It is surprising that
this type of argument works in an infinite-dimensional non-compact
setting, since Moser's trick needs the compactness of the
manifold. Somehow, the compactness of the underlying manifold
makes the job in our case.

This shows that the symplectic geometry of the loop space probably
does not recover the Riemannian geometry of the underlying
manifold in the even dimensional case. This statement will be
clearer after the proof of Theorem \ref{thm:main_unique}, in which
the geometric obstruction for the complete ``uniqueness'' 
of the quasi-symplectic structure is shown.

Finally we discuss how to find contact hypersurfaces in loop (and
path) spaces. The most natural construction is given by

\begin{theorem} \label{propo:expanding}
Assume that the Riemannian manifold $(M,g)$ admits a vector field
$X$ which satisfies $L_X g= g$ and is locally gradient-like, then the lift of $X$ to $\SL(M)$
(respectively $\SP(M)$) is a Liouville vector field for the length
function.
\end{theorem}

This shows in particular that stabilizing the manifold $M$, i.e.
considering $M \times \R$, we obtain contact hypersurfaces in the
loop space.




\section{Symplectic structure} \label{sect:symp}

\subsection{Basic definitions.}
The quasi-symplectic structure in the space of loops of a Riemannian
manifold is defined by taking the exterior differential of the
$1$-form $\mu$ given by equation (\ref{eq:one_form}). To do that we
recall the formula
\begin{equation}
d\alpha (X,Y)= X(\alpha(Y))- Y(\alpha(X)) - \alpha([X, Y]), \label{eq:Cartan}
\end{equation}
which is valid for any $1$-form $\alpha$ and it does not depend on
the vector fields $X,Y$ chosen to extend $X(\gamma)$ and $Y(\gamma)$
for a given point (loop) $\gamma\in \SL(M)$. In our case we start
with two vectors $U, V \in \Gamma(S^1,\gamma^* TM) \simeq T_{\gamma}
\SL(M)$. First define
$$ \theta: (-\e, \e)^2 \times S^1 \to M, $$
satisfying:
\begin{enumerate}
\item $\theta(0,0,t)= \gamma(t)$,
\item $\frac{\partial \theta}{\partial u}(0,0,t)= U(t)$,
\item $\frac{\partial \theta}{\partial v}(0,0,t)= V(t)$.
\end{enumerate}
And so define $\gamma'= \frac{\partial \theta}{\partial t}$, $\hat{U}= \frac{\partial \theta}{\partial u}$ and
$\hat{V}=\frac{\partial \theta}{\partial v}$. They clearly satisfy
\begin{equation}
[\hat{U}, \hat{V}] =0, \label{eq:commuting_exten}
\end{equation}
since they are derivatives of the coordinates of a parametrization. This allows us to compute
\begin{eqnarray}
\hat{U}(\mu(\hat{V}))& = & \frac{d}{du}\left( \frac12 \int_0^1 g_{\theta(u,0,t)}\Big(
\frac{\partial \theta}{\partial v}(u,0,t), \frac{\partial \theta}{\partial t}(u,0,t)\Big) dt \right) \nonumber \\
& = & \frac12 \int_0^1 (g_{\theta(0,0,t)} ( \nabla_{\hat{U}}
\hat{V}, \gamma'(t))+ g_{\theta(0,0,t)} ( V, \nabla_{\hat{U}} {\gamma'}))dt. \label{eq:half_comp}
\end{eqnarray}
In the same way, we obtain
\begin{equation*}
\hat{V}(\mu(\hat{U})) = \frac12 \int_0^1 (g_{\theta(0,0,t)} (
\nabla_{\hat{V}} \hat{U}, \gamma'(t))+ g_{\theta(0,0,t)} ( U,
\nabla_{\hat{V}} {\gamma'}))dt.
\end{equation*}
We are using the torsion-free Levi-Civita connection for the
computations, so $\nabla_{\gamma'} \hat{U}=\nabla_{\hat{U}}\gamma'$
and $\nabla_{\gamma'} \hat{V}=\nabla_{\hat{V}}\gamma'$. Also
$\nabla_{\hat U} \hat{V}- \nabla_{\hat V} \hat{U}= [\hat{U},
\hat{V}] =0$. We shall use the notation $\nabla_{\gamma'} U=
\frac{\partial U}{\partial t}$. So we have, by applying the formula
(\ref{eq:Cartan}), that
\begin{equation}\label{eq:2form_pre}
\begin{aligned}
\omega(U,V) & =  \omega(\hat{U}, \hat{V}) = d\mu(\hat{U}, \hat{V}) =
\hat{U}(\mu(\hat{V}))- \hat{V}(\mu({\hat{U}})) = \\
& =   \frac12 \int_0^1 \Big(
 g_{\gamma(t)}\big(V, \frac{\partial U}{\partial
t}\big) - g_{\gamma(t)}\big(U, \frac{\partial V}{\partial t}\big)
\Big)dt.
\end{aligned}
\end{equation}
Moreover we have
\begin{equation}
0 = \int_0^1 \left( \frac{d}{dt} g(U, V) \right) dt = \int_0^1
\Big(g\big(\frac{\partial U}{\partial t}, V\big)+ g\big(U,
\frac{\partial V}{\partial t}\big)\Big)dt, \label{eq:Stokes}
\end{equation}
which implies
\begin{equation}
\omega(U, V)= \int_0^1 g\big(\frac{\partial U}{\partial t},
V\big)dt. \label{eq:2form}
\end{equation}
Now the kernel of this $2$-form at a point $\gamma$ is given by the
parallel vector fields along $\gamma$. Therefore $\dim \ker(\gamma)
\leq n$. There are several ways of removing the kernel. The simplest
one is to fix a point $p\in M$ and to define
$$ \SL_{p}(M)= \{ \gamma \in \SL(M)\ ; \  \gamma(0)=p \}. $$
This forces the tangent vectors to satisfy
$$ X \in T_{\gamma} \SL_p(M) \Rightarrow X \in \Gamma(S^1, \gamma^* TM), \  X(0)=0. $$
Therefore any parallel vector field is null. So the manifold
$\SL_p(M)$ is symplectic.

Extend our space to $\SP(M)$ where it is still possible to repeat
all the previous computations. We highlight the differences. The
equation (\ref{eq:half_comp}) is exactly the same as it is
symmetric. The equation (\ref{eq:2form_pre}) remains also without
changes. We just need to rewrite equation (\ref{eq:Stokes}) which is
not true anymore and so the final expression for the exterior
differential of $\mu$ becomes
\begin{eqnarray*}
d\mu (U, V) & = & \omega(U,V)= \int_0^1 \Big(g\big(\frac{\partial U}{\partial t},V\big) - \frac{1}{2}\frac{d}{dt} g(U,V)\Big) dt \\
& = & \int_0^1 g\big(\frac{\partial U}{\partial t},V\big)dt -
\frac{g(U(1),V(1))-g(U(0),V(0))}{2}.
\end{eqnarray*}
Is is obviously a closed (being exact) form. Let us compute its
kernel. Assume that $X \in \ker (\omega)_\gamma$. Considering
$\omega(X,V)=0$ for all vectors $V\in T_{\gamma} \SP(M)$ with
$V(0)=V(1)=0$, we obtain that
$$
\frac{\partial X}{\partial t}= 0.
$$
Now by choosing all $V\in T_\gamma \SP(M)$ with $V(0) \neq 0$ and
$V(1)=0$, we conclude that $X(0)=0$. By parallel transport, $X=0$
and so the kernel of $\omega$ is trivial. Hence this form is
symplectic. This proves Proposition \ref{thm:str_path}.

\subsection{Almost complex structures} \label{subs:almost}

There is a canonical almost-complex structure compatible with
$\omega$ in $\SL_{p} M$. Let us construct it. Given a curve
$\gamma:[0,1] \to M$, denote $P_s^t$ the parallel transport
isometry along $\gamma$. There is an isometric isomorphism between
$\gamma^* TM$ and the trivial $T_{\gamma(0)}M$ bundle over $I$
with constant metric $g_{\gamma(0)}$. This allows to translate any
section $U(t) \in \gamma^* TM$ to a section
$P_t^0(U(t))=\hat{U}(t) \in T_{\gamma(0)}M$. This gives rise to a
``d\'eveloppement'' map
  $$
  T_\gamma \SL_p \cong \SL_0 (T_{\gamma(0)} M) .
  $$
Note that if we apply this map to $\gamma'(t)$ itself, we get
$x(t)\in T_{\gamma(0)}M$. Now we define
$$ a(t)= \int_0^t x(s)ds, $$
which is known as the ``d\'eveloppement de Cartan'' of the curve
$\gamma$ in the tangent space $T_{\gamma(0)}M$.  As the covariant
derivative along $\gamma$ becomes the ordinary derivative in
$T_{\gamma(0)}M$, we have that $\gamma$ is a geodesic just when
its d\'eveloppement de Cartan is a line.

Define the almost complex structure $\hat{J}$ in
$T_{\gamma}\SL_p(M)$ as follows: take any vector field $U \in
\Gamma(S^1, \gamma^* TM)$ and compute its ``d\'eveloppement''
$\hat{U}$. Recall that $\hat{U}(0)=\hat{U}(1)=0$ since $U(t)\in
T_\gamma\SL_p(M)$. Fixing an isomorphism $T_{\gamma(0)} M\cong\R^n$,
we have $\hat{U}(t)\in \SC^{\infty}(S^1, \R^n)$. Take its Fourier
series expansion,
$$ \hat{U}(t) = \sum_{k=-\infty}^{\infty} a_k e^{2\pi ikt}, $$
where $a_k\in \C^n$ and $a_{-k} = \bar{a}_k$. Then define
\begin{equation}
\tilde{J}(\hat{U})(t)= \sum_{k<0} (-ia_k)\, e^{2\pi ikt} +
a_0+\sum_{k>0} ia_k\, e^{2\pi ikt}, \label{eq:almost_comp}
\end{equation}
We substract the constant vector $J(\hat{U})(0)$ to get the
almost-complex structure. So we have
$$ \hat{J} (\hat{U})(t) = \tilde{J}(\hat U)(t)- \tilde{J}(\hat U)(0) \in T_\gamma \SL_p(M). $$
To check that it is an almost complex structure we compute
\begin{eqnarray*}
\hat{J}\hat{J}(\hat U) & = & \hat{J} (\tilde{J}(\hat U)- \tilde{J}(\hat U)(0))= \\
 & = & \tilde{J}\tilde{J}(\hat U)- \tilde{J}(\hat U)(0) - \tilde{J}\tilde{J}(\hat U)(0) + \tilde{J}(\hat U)(0)= \\
 & = & \tilde{J}\tilde{J}(\hat U) - \tilde{J}\tilde{J}(\hat U)(0) = \\
 &=&  (-\hat U +2a_0) - (-\hat{U}(0) +2a_0) = -\hat{U}.
\end{eqnarray*}
So we obtain an almost complex structure on $\SL_0 (T_{\gamma(0)}
M)$, then by using the ``d\'eveloppement'' we have an almost
complex structure in $T_\gamma \SL_p(M)$, that is, in $\SL_p(M)$.

To check that $\hat J$ is compatible with the symplectic form
$\omega$ we just compute the value of $\omega$ when trivialized in
the ``d\'eveloppement'', to obtain
\begin{equation}
\omega(\hat{U}, \hat{V}) = \omega \left( \sum_p a_p e^{2\pi ipt}, \sum_q
b_q e^{2\pi iqt} \right) = \sum_k 2\pi i k \Re \langle a_k, {b}_k
\rangle, \label{eq:hat_symp}
\end{equation}
where $a_k, b_k \in \C^n$, and $\langle \, , \, \rangle$ is the
standard Hermitian product in $\C^n$. The associated metric
 $$
 g(\hat{U}, \hat{V})= \omega(\hat{U}, J\hat{V})=
\sum_{k>0} 2\pi k \Re (\langle a_k, {b}_k \rangle + \langle
a_{-k}, {b}_{-k} \rangle)
 $$
is clearly Riemannian. Moreover, $\hat{J}$ is smooth. To check it, recall that
the map
\begin{eqnarray*}
\SF: \{f\in \SC^{\infty}(S^1, \R^n)\ ; \ f(0)=0 \} & \to & \SS \subset (\C^n)^{\infty} \\
f & \mapsto & (a_1, a_2, a_3, \ldots),
\end{eqnarray*}
where $\SS$ is the Schwartz space of sequences of vectors in $\C^n$
with decay faster than polynomial, and $\{ a_k \}$ are the Fourier
coefficients of $f$, is a topological isomorphism (we take in $\SS$
the Fr\'echet structure given by the norms $||(a_k)||_t=\sum k^t
|a_k|$). (Note that the Fourier coefficients satisfy
$a_{-k}=\bar{a}_k$ and $a_0=-\sum_{k\neq 0} a_k$.) The map $\hat{J}$
is conjugated under $\SF$ to the map
\begin{eqnarray*}
\SJ: \SS & \to & \SS, \\
(a_1, a_2, \ldots) & \mapsto & (ia_1, ia_2, \ldots),
\end{eqnarray*}
which is smooth (actually an isometry). So $\hat{J}$ is smooth. This
corrects the folklore statement saying that this almost complex
structure is not smooth in general (see \cite[pag.\ 355]{Wu95}).
\medskip

In the case in which $M$ is a Lie group, there is an alternative
way of defining an almost complex structure for the space of loops
based at the neutral element $e\in G$. To do it we just use the
left multiplication to take the tangent space $T_{\gamma}\SL_e(G)$
to $\Gamma(S^1, T_eG)$, so we obtain an isomorphism $T_{\gamma}
\SL_e(G) \simeq \Gamma(S^1, \R^n)/{\R^n}$, preserving the metric
by construction (the quotient is by the constant maps). So every
particular vector field $X\in T_{\gamma} \SL_e(G)$ is transformed
via the isomorphism to a loop in $\R^n$. Recall that the
isomorphism does not coincide with the one induced by the
``d\'eveloppement''  unless the group is flat (an
abelian group). Once we have set up the previous identification,
the formula (\ref{eq:almost_comp}) provides again an almost
complex structure. Again we remark that it does not coincide with
the previous one in the cases when both are defined.

\subsection{Uniqueness}
One may ask how canonical the symplectic structures on the loop
spaces are. Such symplectic structure $\omega$ is associated to a
metric $g$ on $M$. Recall that the space of Riemannian metrics on
a manifold is connected. Therefore, the Morse trick may help to
prove that the associated loop spaces are all symplectomorphic.
Recall that the Moser's trick for exact symplectic forms works as
follows. Assume that the $1$-parametric family of forms $\omega_t=
d\mu_t$ are symplectic on a manifold $N$. We want to find $\phi_t:
N \to N$, such that $\phi_t^* \omega_t= \omega_0$. Let $Y_t$ be
the vector field generating $\phi_t$. Then
 $$ L_{Y_t} \omega_t = d \big(\frac{d\mu_t}{dt}\big) $$
is satisfied. This is true if
\begin{equation}
i_{Y_t} \omega_t = \frac{d\mu}{dt}, \label{eq:Moser}
\end{equation}
which has a unique solution since $\omega_t$ non-degenerate. In
our case we are given two different Riemannian metrics $g_0$ and
$g_1$ in $M$, so there is a path of metrics $g_t$ that joins them.
Therefore there is a path of exact symplectic forms $\omega_t$ in
$\SL_p(M)$. Now we substitute into equation (\ref{eq:Moser}) in
our case to obtain
\begin{equation}
\int_0^1 g_t\Big(\frac{\nabla Y_t(\gamma)(s)}{ds}, X(s)\Big)ds =
\frac{d}{dt} \int_0^1 g_t(X(s), \gamma'(s)) ds,
\end{equation}
for all $X\in \Gamma(S^1, \gamma^* TM)$. (Note that $Y_t(\gamma)\in
T_{\gamma}\SL_p(M)$.) This equation does not have continuous
solutions in general. This is because we can compute $\frac{\nabla
Y_t(\gamma)(s)}{ds}$ in a continuous way and we may assume that
$Y_t(\gamma)(0)=0$, but then there is no reason to expect that
$Y_t(\gamma)(1)=0$. So, it turns out that we do not get a uniqueness
result for the symplectic structure. We will see a way of partially
avoiding this obstruction. For this we are forced to change our
point of view and work over the space $\SL(M)$, where the forms
$\omega$ are quasi-symplectic.

\bigskip

{\bf Proof of Theorem \ref{thm:main_unique}.} \newline We try to apply Moser's trick in the space $\SL(M)$. Recall that
$\omega$ is not symplectic in this case, since it has a finite
dimensional kernel. Denote by $\SG(M)$ the space of Riemannian
metrics over $M$.

Given a point $p \in M$, define
$$ \so_0(p) = \{ a \in \so(T_pM)\ ; \ \det(a)=0 \}.$$
Recall that $\so_0(p)\neq \so(T_pM)$ if $\dim(M)$ is even. In that
case it is a codimension $1$ (singular) submanifold. This defines a fibration
$$ \so_0(M) \to M, $$
which is a (non-linear) subbundle of the bundle $\so(TM)$.

For a metric $g \in \SG(M)$, the curvature $R_g$ associated to the
metric is a section of the bundle $\so(TM)\otimes \Omega^2(M)$,
therefore it defines a bundle map
 $$
 R_g: \bigwedge{}^2(TM)\to \so(TM)\, .
 $$
Note that these bundles have the same rank. We say that a metric $g$
is R-generic if for each $p\in M$ there are two vectors $u,v\in
T_pM$ such that $\det(R_g(u,v))\neq 0$, in the case $\dim M\geq 4$,
or if $R_g$ vanishes only at a discrete set of points, in the case
$\dim M=2$.


Define
 $$
 \SG_0(M)= \{ g \in \SG(M): g \, \text{ is R-generic} \}.
 $$
Given two metrics $g_0'$ and $g_1'$, there are two metrics $g_0$ and
$g_1$ such that $g_i$ is as close as desired to $g_i'$ and $g_i$ is
R-generic. Moreover we may do the same for paths $g_t$ of R-generic
metrics. This is clear in the case $\dim M=2$. If $n=\dim M\geq 4$,
we work as follows: for each $g\in \SG$ and $p\in M$, there exists a
small perturbation of $g$ which is R-generic in a neighborhood of
$p$. This is true since the subspace
$\Hom(\bigwedge^2(T_pM),\so_0(p))\subset \Hom (\bigwedge^2
(T_pM),\so(T_pM))$ is of codimension ${n\choose 2}$, which is bigger
than $n$. The result follows from a Sard-Smale lemma applied to the
functional
\begin{eqnarray*}
\SG(M) \times M & \to & \Hom (\bigwedge{}^2 (T M),\so(TM)) \\
(g, p) & \to & R_g(p).
\end{eqnarray*}

Now define the set
$$ \SO_0(T_pM)= \{ A \in \SO(T_pM)\ ; \ A-\Id \text{ is not invertible} \}. $$
This is a codimension $1$ stratified submanifold of $\SO(T_pM)$ and
defines a bundle
$$ \SO_0(TM) \to M. $$
Let us define
$$ \SL_s (M)= \{ \gamma \in \SL(M)\ ; \ P_0^1- \Id \ \text{ is not invertible} \}. $$
There is a map
\begin{eqnarray*}
\theta: \SL(M) & \to & \SO(TM) \\
\gamma & \to & (P_0^1)_{\gamma} \in \SO(T_{\gamma(0)}M).
\end{eqnarray*}
Clearly $\SL_s(M) = \theta^{-1} (\SO_0(TM))$. Therefore if $\theta$
is generic in a suitable sense, $\SL_s(M)$ will be a stratified
codimension $1$ submanifold. We claim that if $g \in \SG_0(M)$, this
is the case. To check this, pick $\gamma \in \SL_s(M)$. Being $g$ an
R-generic metric at $T_{\gamma(0)}M$, there exist two vectors
$u,v\in T_{\gamma(0)}M$ such that $\det(R_g{(u,v)}) \neq 0$.

Recall \cite[Subsection 15.4.1]{Be02} that if we extend $u,v$ to a
neighborhood of $\gamma(0)$ in such a way that they define a local
pair of coordinates $(x,y)$ where
$$\frac{\partial}{\partial x}\big|_{(0,0)}=u,$$
$$\frac{\partial}{\partial y}\big|_{(0,0)}=v,$$
and we define  the path $\gamma_{u,v}^s$ as (in the coordinates $(x,y)$):
\begin{itemize}
\item $\gamma_{u,v}^s(t)= ( 4st, 0)$, $t\in [0, 1/4]$,
\item $\gamma_{u,v}^s(t)= ( s, 4s(t-1/4))$, $t\in [1/4, 1/2]$,
\item $\gamma_{u,v}^s(t)= ( 4s(3/4-t), s)$, $t\in[1/2, 3/4]$,
\item $\gamma_{u,v}^s(t)= ( 0, 4s(1-t))$, $t\in [3/4, 1]$,
\end{itemize}
we obtain
$$ \lim_{s \to 0} \frac{(P_0^1)_{\gamma_{u,v}^s}}{s} = R_g{(u,v)}. $$
Take the path which is the juxtaposition of $\gamma$ with
$\gamma_{u,v}^s$,
$$ \beta_s = \gamma \ast \gamma_{u,v}^s. $$
This family of paths determines a tangent vector in $T_{\gamma} \SL(M)$. We will show
that it is transverse to the submanifold $\SL_s(M)$. The holonomy of
$\beta_s$ is
\begin{equation}
(P_0^1)_{\beta_s}= P_0^1\left( \Id+ s\,R_g({u,v}) \right)+O(s^2),
\label{eq:defmono}
\end{equation}
where $P_0^1$ is the holonomy of $\gamma$. Now embed $\gl(n,\R)
\subset \gl(n,\C)$ by complexifying. Then $P_{0}^{1} \in \gl(n,\C)$
admits a Jordan canonical form
$$ J=B P_0^{1} B^{-1}, $$
where $B \in \GL(n, \C)$. Multiply on the left and on the right by
$B$ and $B^{-1}$ the expresion (\ref{eq:defmono}) to obtain
\begin{equation}
B (P_0^1)_{\beta_s} B^{-1}= J + s\, B P_0^1 R_g{(u,v)}
B^{-1}+O(s^2). \label{eq:jordan}
\end{equation}
Now, it is easy to check that the eigenvalues of the right hand side
of (\ref{eq:jordan}), for $s$ small enough, are far away from zero
or grow faster than $\epsilon s$, for some fixed $\epsilon>0$. Since
$B$ is just a change of coordinates matrix for $(P_0^1)_{\beta_s}$
on the left hand side of (\ref{eq:jordan}), the eigenvalues of
$(P_0^1)_{\beta_s}$ are the same than those of the right hand side
of (\ref{eq:jordan}). This implies that the tangent vector
determined by $\beta_s$ is transverse to $\SL_s(M)$.

In the case $\dim M=2$, we work as before when $\gamma(0)$ does not
coincide with either of the points $p\in M$ with $R_g(p)=0$. We can
define $\SL_s(M)$ to be the same set as before together with the
paths starting at a point $p$ with $R_g(p)=0$. This latter set is of
codimension $n$, so $\SL_s(M)$ is still of codimension $1$, as needed, in this
case.

Recall that it is fundamental for this argument to work that the
dimension of $M$ is even. If the dimension is odd, then $\SL_s(M) =
\SL(M)$.

\medskip

All the computations done at the beginning of this subsection remain
valid for the family $\{ g_t \}$ and so we can follow Moser's trick
to end up with equation (\ref{eq:Moser}) again. Now we recall that
the solution $Y_t$ is not unique, basically because we are working
with $\SL(M)$, and there is a space of parallel vector fields along
$\gamma$, which are in the kernel of the quasi-symplectic form
$\omega$. In general we obtain $Y_t^w=Y_t+ w$ as a valid solution
where $Y_t$ is a particular solution and $w$ is a parallel vector
field along $\gamma$ (here we may have $w(0) \neq w(1)$ -- actually,
this will be the case). Again in general $Y_t^w(1) \neq Y_t^w(0)$.
However a careful choice of $w$ may help. We fix the equation
$$ Y_t(0)+w(0)=Y_t(1)+w(1), $$
that if solved for some $w(0) \in \R^n$, leads to a
smooth solution of equation (\ref{eq:Moser}). The previous equation leads
to
\begin{equation}
P_0^1(w(0))-w(0)= Y_t(1)- Y_t(0), \label{eq:adjust_cont}
\end{equation}
that clearly has a unique solution whenever $\gamma \not \in
\SL_s(M)$. So equation (\ref{eq:Moser}) has a unique continuous
solution outside a set of positive codimension in $\SL(M)$. We are
aiming to construct an ``approximate'' solution to
(\ref{eq:adjust_cont}). To get this, we can perturb equation
(\ref{eq:adjust_cont}) to
\begin{equation}
(\lambda P_0^1-\Id)(w_{\lambda}(0))= Y_t(1)- Y_t(0),
\label{eq:adjust_cont2}
\end{equation}
which always admits a solution for $|\lambda|<1$ since $P_0^1 \in
\text{SO}(n)$. Now we assume that $\lambda$ is a smooth map
$\lambda: \SL(M) \times [0,1] \to [1-\epsilon,1]$ satisfying
\begin{enumerate}
\item $\lambda (\gamma, t)= 1-\epsilon$ if $\gamma \in \SL_s(M)$
for the metric $g_t$. \item $\lambda (\gamma, t)= 1$ on a small
neighborhood of $\SL_s(M)$.
\end{enumerate}
This defines a family of flows $\{Y_t^{\lambda} \}$ (for
the given constant $\epsilon>0$). The integrals at time $1$ of
these flows are generating 
the family $\phi_{\epsilon}$ required in the statement of the
theorem.

The flow exists and is unique. This is due to the fact that it can
be understood as a parametric (smoothly dependent) family of flows
in $M$ and there we have existence and uniqueness (for all times).
The smooth dependency in the parameters gives that the flow is
Fr\'echet smooth. We also need to prove that the flow is by
diffeomorphisms since, being $\SC^{\infty}(S^1, M)$ a Fr\'echet
manifold, this is not automatic. But this follows from
$\phi_{-\epsilon} \circ \phi_{\epsilon}=\id$, and so the flow maps
admit inverses.

It is a routine to check that $\phi_{\epsilon}$ takes the
quasi-symplectic form associated to $g_0$ to a form
$\epsilon$-close in $L^2$-norm to the quasi-symplectic form
associated to $g_1$. Moreover, since $g_i$ and $g_i'$ are as close
as needed we can also claim that their associated quasi-symplectic
structures are close to each other in $L^2$-norm. \hfill $\Box$

\medskip

Finally, it is remarkable to note that the possibility of addition
of parallel vector fields has been the key to find a continuously
varying family of functions which are solutions to
(\ref{eq:Moser}) almost everywhere. This is the geometric reason
for which we cannot extend the computation to the case of the
(symplectic) space $\SL_p(M)$.

\section{Loop spaces as contact manifolds.} \label{sect:contact}
We want to check whether the symplectic manifold $\SL_p(M)$ has
hypersurfaces of contact type on it. We prove now Theorem
\ref{propo:expanding}.

\bigskip

{\bf Proof of Theorem \ref{propo:expanding}.} \newline Let $X$ be
a vector field on $M$ satisfying $L_X g =g$. Then $\nabla X\in
\End(TM)$, and its symmetrization is $\frac12 \Id$. This follows
since, for $Y$, $Z$ vector fields on $M$, we have
 $$
 \begin{aligned}
  g(\nabla_Z X,Y) + &\, g(\nabla_Y X,Z) = \\
  &= g(\nabla_X Z,Y) +g(\nabla_X Y,Z)  - g(L_X Z,Y) -g(L_X Y,Z) =\\
   &= X(g(Y,Z)) + (L_Xg) (Y,Z) - L_X(g(Y,Z)) =\\
   &= g(Y,Z),
  \end{aligned}
  $$
where we have used that $L_XZ=\nabla_XZ-\nabla_ZX$ on the second
line. The anti-symmetrization of $\nabla X$ is $\SA(\nabla
X)=\SA(\nabla X^\#)_\#= (d X^\#)_\#$, where $X^\#$ is the
$1$-form associated to $X$ (``raising the index''), and the
$(\cdot)_\#$ means ``lowering the index'' with the metric. Recall that
a vector field is ''locally gradient-like'' in a neighborhood $U$ for a metric $g$ if it is $g$-dual of some exact $1$-form $df$, where $f$ is a  function $f: U \to \R$. Thus, if $X$
is locally gradient-like, then $X^\#$ is a locally exact, i.e.
closed, $1$-form and so $\SA(\nabla X)=0$. Then $\nabla X=\frac12 \Id$.

Associated to $X$ there is an induced vector field $\hat{X}$ on
$\SL(M)$. It is defined as follows: for $\gamma\in\SL(M)$,
$\hat{X}_\gamma \in T_\gamma \SL(M)$ is given by
$\hat{X}_\gamma(t) =X(\gamma(t))$. We want to check that
$L_{\hat{X}}\mu=\mu$. For $Y\in T_\gamma \SL(M)$, we have
 $$
 \begin{aligned}
  \alpha(Y) &= i_{\hat{X}}\omega(Y)=\omega (\hat{X},Y) =\\
   &= \int_0^1 g\big(\frac{\partial X}{\partial t}, Y\big)
   dt = \int_0^1 g(\nabla_{\gamma'} X, Y)dt=\\
   &= \frac12 \int_0^1 g(\gamma' , Y)dt= \mu (Y) .
  \end{aligned}
  $$
So $\alpha=\mu$. Then
 $$
 L_{\hat X}\mu= di_{\hat X}\mu + i_{\hat X}d\mu =
 di_{\hat X}i_{\hat X}\mu + i_{\hat X}\omega =0+\alpha=\mu.
 $$
From this, it follows that $L_{\hat X}\omega=L_{\hat
X}d\mu=dL_{\hat X}\mu=d\mu=\omega$, as required.\hfill $\Box$

\begin{remark}
The manifolds to which the previous result applies, that is, those
satisfying $L_Xg=g$ with $X$ locally gradient-like, are locally of
the form $(N \times \R, e^t(g + dt^2))$ with expanding vector
field $X=\frac{\partial}{\partial t}$. This follows by writing
$X=\grad f$, with $f>0$ and putting $t=\log(f)$.

Two examples are relevant:
\begin{itemize}
\item $M=N\times \R$, with $(N,g)$ a compact Riemannian manifold.
Give $M$ the metric $e^t(g+dt^2)$. \item Let $(N,g)$ be an open
Riemannian manifold with a diffeomorphism $\varphi:N\to N$ such
that $\varphi^*(g)=e^\lambda\, g$, $\lambda>0$. Then take
$M=(N\times [0,\lambda])/\sim$, where $(x,0)\sim
(\varphi(x),\lambda)$ and $M$ has the metric induced by $e^t(g +
dt^2)$.
\end{itemize}
\end{remark}

To finish, let us check that the familiar finite dimensional picture
translates to this case.

\begin{proposition}
Let $(M,g)$ be a Riemannian manifold which has a locally
gradient-like vector field $X$ satisfying $L_X g = g$. Then the
hypersurface
 $$ \SL_{p,1} (M) = \{ \gamma \in \SL_p (M)\ ; \ \length(\gamma)=1 \} $$
is a contact hypersurface of $\SL_p (M)$.
\end{proposition}

\medskip

{\bf Proof.} \newline We need to check that $\alpha= i_{\hat{X}}
\omega$ is a contact form on $\SL_{p,1} (M)$. First we claim that
$\alpha$ is nowhere zero on that submanifold. If this were not the
case, then we would have that
 \begin{equation}
 (i_{\hat{X}} \omega)|_{\SL_{p,1}(M)} = 0, \label{eq:contact_form}
 \end{equation}
and we know that $i_{\hat{X}} \omega (\hat{X})= 0$. Note that
$\hat{X}$ is transversal to $\SL_{p,1}(M)$ since the flow increases
the length of the loop. Therefore we have that $i_{\hat{X}} \omega =
0$ and $\omega$ would not be symplectic, which is a contradiction.
So $\alpha$ is nowhere zero.

Now we have the distribution $(\ker \alpha, \omega = d \alpha)$ on
$\SL_{p,1}(M)$. To finish we check that it is symplectic. Assume
that $Y \in T_{\gamma} \SL_{p,1}(M)$ satisfies that
 $$ \omega(Y, Z)=0, $$
for all $Z \in T_{\gamma} \SL_{p,1}(M)$. Moreover we have that
$\omega(Y, \hat{X})= -\alpha(Y)=0$. Hence $i_Y \omega=0$, and we get
a contradiction.


\begin{corollary}
Given a Riemannian manifold $(M, g)$, the manifold $(M\times \R,
e^{\lambda}(g+ d\lambda^2))$ has an associated space of loops of
length one with a canonical contact form. For a loop $\gamma$ and
$Y$ vector field along $\gamma$, we denote $\gamma= (\gamma_1,
\gamma_2)$ and $Y=(Y_1,Y_2)$ according to the decomposition $M
\times \R$. Then the contact form is given by
 $$
 \alpha(Y) =\mu(Y) =\frac12
 \int_0^1 e^{\gamma_2(t)} \left( g(\gamma_1'(t),Y_1(t))+
 \gamma_2'(t)Y_2(t) \right) dt. $$
\end{corollary}


\subsection{Reeb vector fields.}
We compute the Reeb vector field associated to the contact form. We
need to do it in $\SL_{1}(M)$, that is a quasi-contact space
(instead of contact). This is necessary in order to obtain a smooth
Reeb vector field.

\begin{lemma}
The Reeb vector field associated to $(\SL_{1} (M), \alpha)$ is the vector
 $$ R = \frac{\gamma'}{||\gamma'||}\cdot h(\gamma), $$
where $h: \SL_{1}(M) \to \R$ is a $\Diff(S^1)$-invariant strictly
positive function on the loop space.
\end{lemma}

{\bf Proof.} \newline The condition for a vector $V\in T_\gamma \SL
(M)$ to belong to $T_\gamma \SL_{1}(M)$ is
 $$ \lim_{s\to 0} \frac{\int_0^1 || \gamma' + s\frac{dV}{dt} ||dt}{s} = 0. $$
So we are asking the vector to satisfy
\begin{equation*}
 \int_0^1 g\big(\frac{\partial V}{\partial t}, \gamma'(t)\big) \frac1{||\gamma'||} dt =0.
\end{equation*}
Thus the previous equation can be rewritten as
 $$ \omega\Big(V, \frac{\gamma'}{||\gamma'||}\Big)=0. $$
Imposing this condition for every non-zero vector $V \in
T_\gamma\SL_{1}(M)$, we get that the Reeb vector field is a positive
multiple of $\frac{\gamma'}{||\gamma'||}$, so proving the statement.
It remains to be checked that $h(\gamma)$ is invariant by changes of
parametrization preserving the origin. This is a consequence of the
invariance under changes of coordinates of the defining equation
 $$
 \alpha (R)= i_{\hat{X}} \omega \Big( \frac{\gamma'}{||\gamma'||}
 h(\gamma)\Big) = 1.
 $$
\hfill $\Box$

There are more solutions to the Reeb vector field equation since for
any parallel vector field $w$ along $\gamma$, we can add it to $R$,
so that $R+w$ defines another solution to the equation. However
those solutions are not continuous (as functionals on $\SL(M)$) in
general. This is the reason that prevents us to define the Reeb
field in the space $\SL_{p,1}(M)$. Again the flexibility of the
quasi-contact structure is the key to find the solution.

We have the following
\begin{lemma}
All the Reeb orbits of $\alpha$ are closed. For a loop
\begin{eqnarray*}
\gamma: S^1 & \to & M \\
t & \to & \gamma(t),
\end{eqnarray*}
the Reeb orbit passing through it has period equal to the length of
$\gamma$ divided by $h(\gamma)$.
\end{lemma}

{\bf Proof.} \newline Take $\gamma:S^1 \to M$. The arc-length
parametrization of $\gamma(S^1)$ is denoted as $\gamma_p$, where $p$
is the point of $\gamma(S^1)$ in which the arc-length parameter
starts. So we define $\theta(s,t)= \gamma_{\gamma(t)} (h(\gamma) s)$
which is clearly the Reeb orbit starting at $\gamma$. It is periodic
with period $\length(\gamma)/ h(\gamma)$. \hfill $\Box$

\end{document}